\documentclass[11pt]{amsart}
\usepackage{amsmath}
\usepackage{amssymb}

\newtheorem{proposition}{Proposition}
\newtheorem{lemma}{Lemma}

\begin{document}
\title[Zeta representations of groups]
{Generalized  Zeta function representation of groups and 2-dimensional Topological Yang-Mills theory:\\
The example of $GL(2,{ \mathbb F}_q)$ and $PGL(2,{ \mathbb F}_q)$.}
\author{Ph.Roche}
\address{Ph.Roche, Universit\'e Montpellier 2, CNRS,  L2C, IMAG}
\email{philippe.roche@univ-montp2.fr}
\begin{abstract}
We recall the relation between Zeta function representation of groups and two-dimensional topological Yang-Mills theory through Mednikh formula.
We prove various generalisations of Mednikh formulas and  define generalization 
of Zeta functions representations  of groups. We compute  some of these functions in the case of the finite group $GL(2, {\mathbb F}_q)$ and $PGL(2,{\mathbb F}_q).$
We recall the table characters of these groups  for any $q$,  compute the Frobenius-Schur indicator of their irreducible representations and  give   the explicit structure of their  fusion rings. \end{abstract}

\maketitle

\section{Introduction}
$2$-dimensional  Yang-Mills theory associated to compact real Lie groups is now a well studied subject initially  developped  in  \cite{Mi}\cite{Wi}. In the case where the gauge group  $G$ is a connected simply connected compact Lie group (for example $SU(n)$) , the partition function of $2$-dimensional Yang Mills theory on a Riemann surface $\Sigma$ of genus $g$ can be computed explicitely and is equal to
\begin{equation}
{ \mathcal Z}_G(\lambda,g)=\sum_{\pi\in \widehat{G}} \frac{1}{(dim \pi)^{2g-2}}\exp(-\lambda a(\Sigma) C(\pi))
\end{equation}
 where the summation is over the set of  classes of  irreducible finite dimensional complex representations of $G,$ 
$\lambda$ denotes the Yang-Mills coupling constant, $a(\Sigma)$ the area of $\Sigma$ and $C(\pi)$ the value of  the Casimir element of $\pi$.

In the limit where $\lambda$ goes to zero, the so-called topological limit, the partition function depends only on the genus and is equal to:
\begin{equation}\label{partitionfunctiontopo}
{ \mathcal Z}_G(g)=\sum_{\pi\in \widehat{G}} \frac{1}{(dim \pi)^{2g-2}}.
\end{equation}

The initial motivation of the present work was to study low dimensional field theory in the case where the group is a $p-$adic one and to study generalisation of the previous formulas in the context of $p-$adic compact Lie group (such as for example  $SL(n, \mathbb{Z}_p)).$ This could be the initial step to build adelic low dimensional field theory, theories which could be of fundamental importance. See  Y.Manin \cite{Ma} for reflections on this topic and    \cite{BrFr} for a review of $p$-adic and adelic physics.

The formula (\ref{partitionfunctiontopo})  gives a bridge with the study of the behaviour of dimensions of representations of groups.
Let $G$ be a compact topological group, if $\pi$ is  a continuous complex representation of necessarily finite dimension $dim \pi$, we denote  $\chi_\pi$ its character and $V_\pi$ the vector space on which $\pi$ acts. 

Let $\widehat{G}$ be the set of isomorphism 
 classes of irreducible finite dimensional complex representations. We assume that for each positive 
integer $n$ the number of elements of $\widehat{G}$ having dimension $n$ is finite of cardinal
 $a_n.$ Groups having these properties are called rigid \cite{Voll1}.

One can define the $\zeta_G$ function representation  of $G$ by
\begin{equation}
\zeta_G(s)=\sum_{n\geq 1}\frac{a_n}{n^s}=\sum_{\pi\in \widehat{G}} \frac{1}{(dim \pi)^{s}}, \,\, s\in \mathbb{C}.
\end{equation}

In the case where $G$ is a connected simply connected real  compact Lie group, $G$ is rigid and  the function $\zeta_G$ is completely determined by the root system of $G$. Fix a Cartan torus $T$  of $G$, it defines a root system $\Phi$, fix a choice $\Phi^+$ of polarization of $\Phi$ , there is a one to one correspondence between $\widehat{G}$ and the set  of positive integer weights $\Lambda^+.$  If $\lambda\in\Lambda^+, $  it defines an irreducible representation $\pi_\lambda$ of $G$ which dimension is given by the formula $dim(\pi_\lambda)=
\prod_{\alpha\in \Phi^+}\frac{\langle\alpha,\lambda+\rho\rangle}{\langle\alpha,\rho\rangle }$ where $\rho$ is  half of the sum of positive roots.
The $\zeta_G$ representation function has therefore the  simple expression:
\begin{equation}
\zeta_G(s)=\sum_{\lambda\in \Lambda^+}\prod_{\alpha\in \Phi^+} \left(\frac{\langle\alpha,\rho\rangle }{\langle\alpha,\lambda+\rho\rangle}\right)^s, \,\, s\in \mathbb{C}.
\end{equation}
It has been shown in \cite{LL} that the abscisse  of convergence  of $\zeta_G$, i.e the infimum of $\alpha>0$ for which $\zeta_G(s)$ is absolutely convergent for $Re(s)>\alpha$, is given by  $2/h$ where $h$ is the Coxeter number of $G.$

In the case of  $p-$adic compact Lie group the situation is much more complicated. This comes from the fact that the
 classification of finite dimensional irreducible complex representations of $G$ is more involved and uses an analog of Kirillov Theory . For a review of known results  see \cite{Klopsch1, Klopsch2, Voll1}. When $G$ is a compact $p-$adic Lie group which is FAb, i.e $H/[H,H]$ is finite for every open subgroup $H$ of $G$, then $G$ is rigid.
 A fundamental theorem, due to A.Jaikin-Zapirain \cite{JZ}, asserts that, under the assumption that $G$ is a compact p-adic Lie group which is FAb, there exists integers $n_1,...,n_k$ and rational functions $f_1,...,f_k$ such that 
 \begin{equation}
 \zeta_G(s)=\sum_{j=1}^k \frac{f_j(p^{-s})}{n_j^s}.
 \end{equation}
 
 In some very few cases one can even compute exactly $\zeta_G $, this is the case for $SL(2,{\mathcal O})$ \cite{JZ}  where ${\mathcal  O}$ is the ring of integers of a p-adic field and more recently for $SL(3, {\mathcal O })$\cite{AKOV}.

In the first part of this article we will review results of interest on $\zeta_G$ function representations and introduce other functions generalizing these $\zeta_G$ functions which are connected to correlation functions of 2D topological Yang Mills theory field theory and its generalisations.
These functions are still actively studied in the real case but have not received any attention  in the p-adic context.

In the second part of this article we compute explicitely these functions for the simplest case where $G$ is a finite group of type 
$GL(2,{\mathbb F}_q)$ and $PGL(2,{\mathbb F}_q).$ Although these functions are very simple compared to those associated to $p-$adic compact Lie groups,  they are already  non trivial and are the preliminary step   for computing the generalized functions for $p-$adic Lie groups of the type $GL(2, {\mathcal O}),$$  PGL(2, {\mathcal O}), $ where ${\mathcal O}$ is the ring of integer of a $p-$adic field. This is the subject of   a subsequent work \cite{BaRo}.

We also collect in an appendix the character table of the groups
$GL(2,{\mathbb F}_q)$,  $PGL(2,{\mathbb F}_q)$ enriched by the list of Frobenius-Schur indicators as well as the explicit decomposition of tensor product of irreducible representations. This last result, although being straightforwardly computed  from the character table of the groups, has not appeared in such a  detailed form in the litterature (see however \cite{AbPa} \cite{Ka}).

\section{Generalized Zeta function representations of groups}
\subsection{Zeta function representation}

In the case where $G$ is finite, the value of $\zeta_G$ on positive even integer has a nice geometrical 
interpretation. 
We will  use this interpretation to define  new zeta functions which 
will have the property that evaluated on some integers they will also have  nice geometrical 
interpretation.

Let $G$ be a finite group, let $\Sigma$ be a compact connected orientable topological surface of genus $g\geq 0$. We fix a point $x$ on $\Sigma,$
and consider the  set of  group morphisms $X_G(x,\Sigma)=Hom(\pi_1(x,\Sigma),G).$

Because the groups $\pi_1(x,\Sigma)$ do not depend on $x$ up to isomorphism, we will not mention $x$ in a context where  results do not depend on the choice of $x.$

$G$ acts on $X_G(x,\Sigma)$ by conjugation and we denote $X_G(x, \Sigma)/Ad G$ the set of equivalence classes.
$X_G(x,\Sigma)/AdG$ is in bijection with the set of $G$-principal 
bundle over $\Sigma$. This set is also called the character variety of $\Sigma$.

The formula of Mednykh \cite{Med} gives:
\begin{equation}
\frac{\vert X_G(x,\Sigma)\vert }{\vert  G\vert}=\vert G\vert^{-\chi(\Sigma)}\sum_{\pi\in \widehat{G}}
 (dim \pi)^{\chi(\Sigma)},
\end{equation}
where $\chi(\Sigma)=2-2g$ is the Euler characteristic of $\Sigma$.

We recall here a simple proof of this formula.

Proof:

$\Sigma$ is homeomorphic to the connected sum $\Sigma\simeq T^{\# g}$ where $T$ is the two dimensional torus.
Therefore  $\pi_1(x,\Sigma)$ has a presentation by generators $a_1,\cdots,a_g, b_1, \cdots, b_g$  and relation $\prod_{i=1}^g a_i b_i a_i^{-1} b_i^{-1}=e.$

We therefore obtain  $$X_G(x,\Sigma)=\{(A_1,\cdots,A_g,B_1,\cdots, B_g)\in G^{2g}, \prod_{i=1}^g [A_i ,B_i]=e\}.$$
We denote $\theta_{T}$ the function $\theta:G\rightarrow {\mathbb N}, \theta_{T}(k)=\vert \{(A,B)\in G^2, ABA^{-1}B^{-1}=k,\}\vert$.
As a result 
\begin{eqnarray*}
\vert X_G(x,\Sigma)\vert&=&\vert G \vert^{2g}\int_{G^{2g}}\delta( \prod_{i=1}^g [A_i, B_i])\prod_i dA_i \prod_i dB_i\\
&=&\vert G\vert^g\int_{G^{g}}\theta_T (k_1)\cdots\theta_T (k_g)\delta(k_1\cdots k_g)dk_1 ..dk_g\\
&=&\vert G\vert^{g-1} (\underbrace{\theta_T *\cdots *\theta_T}_{g \,\text{times}})(e)\end{eqnarray*}
 where $\int$ is the normalized Haar measure on $G,$ $\delta$ is the delta function at $e$ defined by $\delta(x)=1$ if $x=e$ and zero otherwise, and $*$ 
denotes the convolution product, $(\phi*\psi)(x)=\int_{G}\phi(xy)\psi(y^{-1})dy.$

The function $\theta_T$ admits the following Fourier decomposition:
\begin{equation}
\theta_T=\vert G\vert \sum_{\pi\in \widehat{G} }dim(\pi)^{-1}\chi_\pi,
\end{equation} which follows 
from $\vert G\vert \delta=\sum_{\pi\in \widehat{G} }dim(\pi)\chi_\pi$ and from  the formula 
\begin{equation}\label{orthogonality}
\int_{G}\langle \phi, \pi(g)(v)\rangle \langle \psi, \pi(g^{-1})(w)\rangle dg =\frac{\langle \phi, w\rangle \langle \psi, v\rangle }{dim \pi}, \forall v,w\in V_\pi, \forall \phi,\psi\in V_\pi^*.
\end{equation}
Using the Fourier decomposition of $\theta_T$ and the identity $\chi_\pi*\chi_{\pi'}=\frac{1}{dim \pi}\delta_{\pi,\pi'}\chi_\pi,$ we  obtain Mednykh formula.$\square$

\medskip

Note also that $\zeta_G(0)=\vert  \widehat{G} \vert$ which  is also equal to  the number of conjugacy classes of $G$.  We also have $\zeta_G(-2)=\vert G\vert^2.$

When $G$ is a connected simply connected semisimple compact Lie group, an analogue of this correspondence  between the geometry of $G$-principal bundles  and  zeta function representation exists and was first  obtained by  E.Witten \cite{Wi}:
the evaluation of $\zeta_G$
on the integers $2g-2$ for $g\geq 2$ gives the volume of the moduli space of $G$-flat connections on a compact orientable topological surface of genus $g.$ 
A rigorous  mathematical proof of it was given for example in \cite{MeWo},
one has:
 \begin{equation}
 \label{volumeofflatconnections}Vol( X_G(x,\Sigma)/AdG)=\vert Z(G)\vert (Vol(G))^{2g-2}\zeta_G(2g-2)
 \end{equation}
 where $Z(G)$ is the center of $G$, $Vol(G)$ is the Riemannian volume of $G$ 
computed using the invariant metric obtained from the Killing form.  

Note that it is $X_G(x,\Sigma)/AdG $ which appear and not $X_G(x,\Sigma).$

These results can be generalized to the case of surfaces with boundaries, which is the content of the next section.

\subsection{Generalized zeta functions and surfaces with boundaries}
Let $G$ be a finite group, the $\zeta_G$  function can be generalized to the case where one consider connected compact orientable  topological surfaces with boundaries.
Let $\Sigma$ be a connected compact orientable topological surface $\Sigma$ of genus $g$ and define $\Sigma'= \Sigma \setminus (\cup_{j=1}^r D_j)$ where $D_j$ are  open disks included in $\Sigma$  with $(\overline D_j)$ pairwise non intersecting. We  have $\partial \Sigma'= C_1\cup \cdots \cup C_r,$ where $C_1,..., C_r$ are homeomorphic to circles.  
Let $x$ be an interior point of $\Sigma'$, and fix closed curves $\ell_1,...,\ell_r$ originating from $x$ such that the homotopy class of $\ell_j$ is the same as 
$C_j$.
The homotopy group  $\pi(\Sigma',x)$ is isomorphic to the group generated by $a_i,b_i,m_j, i=1,\cdots  g,  j=1,\cdots r,$ with relation
 $ \prod_{i=1}^g a_i b_i a_i^{-1} b_i^{-1}\prod_{j=1}^r m_j =e.$

Let denote $\gamma_1$,...,$\gamma_r$ elements of $G$, the  set  of principal $G$-bundles over $\Sigma'$  such that the holonomy  along $\ell_j$  belongs to the conjugacy class of $\gamma_j,  j=1,\cdots r, $ is a finite set denoted  ${X}_G(x,\Sigma,\gamma_1,...,\gamma_r).$ 
If $\gamma\in G$ we denote $O(\gamma)$ the conjugacy class of $\gamma.$ The characteristic function of $O(\gamma)$ is denoted $ \delta_{O(\gamma)}$ and admits the expansion:
\begin{equation}
\label{fouriercaracconjugacyclass}\frac{\vert G\vert}{\vert O(\gamma)\vert }\delta_{O(\gamma)}=\sum_{\pi\in \widehat{G} }\chi_\pi(\gamma^{-1}) \chi_\pi.
\end{equation}

Mednykh formula can be generalized as:

\begin{equation}
\label{medformulainsertion}\frac{\vert {X}_G(x,\Sigma,\gamma_1, ...,\gamma_r)\vert }{\vert  G\vert}=\vert G\vert^{2g-2}\prod_{j=1}^r\vert O(\gamma_j)\vert
\sum_{\pi\in \widehat{G} }\frac{\prod_{j=1}^r \chi_{\pi}(\gamma_j)}{(dim \pi)^{2g-2+r}}.
\end{equation}

Proof:

It  is a consequence of
\begin{eqnarray*}
&&{X}_G(x,\Sigma,\gamma_1,...,\gamma_r)=\\
&&\{(A_1,..,A_g,B_1,.., B_g,M_1,.., M_r)\in G^{2g+r}, \prod_{i=1}^g [A_i, B_i ]\prod_{j=1}^r M_j=e,\forall j=1,..,r,
M_j\in O(\gamma_j) \},
\end{eqnarray*}
hence
\begin{eqnarray*}
&&\vert {X}_G(x,\Sigma,\gamma_1,...,\gamma_r)\vert=\\
&&\vert G \vert^{2g+r}\int_{G^{2g+r}}\delta( \prod_{i=1}^g [A_i, B_i]\prod_{j=1}^r M_j)\prod_{j=1}^r\delta_{O(\gamma_j)} (M_j) \prod_i dA_i \prod_i dB_i \prod_j dM_j.
\end{eqnarray*}
The formula (\ref{medformulainsertion}) is obtained by using the expansion (\ref{fouriercaracconjugacyclass})  of the characteristic function of conjugacy class and using repetitively the formula (\ref{orthogonality}).
$\square$

\medskip
Therefore we are led to define, for a  rigid compact  group $G$,  a  generalized zeta functions $\zeta_G^{(r)}:G^r\times
 {\mathbb C}\rightarrow {\mathbb C}$ defined by 
\begin{equation}\zeta_G^{(r)}(\gamma_1, ...,\gamma_r) (s)=
\sum_{\pi\in \widehat{G} }\frac{\prod_{j=1}^r \chi_{\pi}(\gamma_j)}{(dim \pi)^{s+r}}, s\in {\mathbb C}, \gamma_1,...,\gamma_r\in G.
\end{equation}
We will call the function $\zeta_G^{(r)}$, the generalized $\zeta$ function representation of $G$ with r-insertions.

In the case where $G$ is a compact connected simply connected real Lie group,  the study of these $\zeta$ function with r-insertions,  has led to numerous work.
First of all the values of these series when $s$ is an even positive integer are connected to the volume of the moduli space of flat connection on a surface with boundaries  with prescribed 
valued of the holononomy in a conjugacy class.  See \cite{Wi} for the statement of the theorem and \cite{MeWo} for a rigorous  proof of it.
When $s$ is an even positive integer, the study of the function $(\gamma_1,...,\gamma_r)\mapsto\zeta^{(r)}(\gamma_1,...,\gamma_r)(s)$ is developped in 
\cite{Sz}, where the   relation of these functions with generalization of Bernouilli polynomials is studied. For recent works on this subject see  \cite{KMT}\cite{BaBoVe}.
In particular it is shown that if $s\in 2{\mathbb N}$  and if we write $\gamma_j=exp(v_j)$ with $v_j\in {\mathfrak h}$ where ${\mathfrak h}$ is a fixed Cartan Lie    subalgebra   of ${\mathfrak g }$, the function
 ${\mathfrak h}^r\rightarrow {\mathbb C}, (v_1,\cdots,v_r)\mapsto \zeta^{(r)}(e^{ v_1},\cdots, e^{v_r})(s)$ is a "locally" polynomial function, meaning that there is a family of open sets
 $(U_i)$ of ${\mathfrak h}^r$  such that $\cup_{i\in I} \overline{U_i}={\mathfrak h}^r,$ and the restriction of   $\zeta^{(r)}(., s)$ to $U_i$  is a polynomial function.

This is easily seen in the case where  $G=SU(2)$: we denote $\pi_n$ the representation of $G$ of dimension $n+1, n\in {\mathbb N}.$ We denote $T$ the Cartan torus of diagonal matrices and let  $t=\begin{pmatrix}1&0\\0&-1
\end{pmatrix} $, we have $\chi_{\pi_n}(e^{i\theta t})=\frac{sin((n+1)\theta)}{sin(\theta)}.$
Therefore, for $s$ non negative  even integer, we have 
\begin{eqnarray*}
  &&\zeta^{(r)}(e^{i\theta_1 t},\cdots, e^{i\theta_r t})(s)=\sum_{n\in {\mathbb N}}\frac{1}{(n+1)^{s+r}}\prod_{j=1}^{r}\frac{sin((n+1)\theta_j)}{sin(\theta_j)}\\
&=&\sum_{m\in {\mathbb N}\setminus\{ 0\}}\frac{1}{m^{s+r}}\prod_{j=1}^{r}\frac{e^{im\theta_j}-e^{-im\theta_j}}{2isin(\theta_j)}\\
&=&\frac{1}{(2i)^r\prod_{j=1}^{r}sin(\theta_j)}\sum_{\epsilon_1,..,\epsilon_r\in\{\pm 1 \}}   \prod_{j=1}^r\epsilon_j \sum_{m\in {\mathbb N}\setminus\{ 0\}}\frac{1}{m^{s+r}}\exp(im\sum_{j=1}^r\epsilon_j\theta_j)\\
&=&\frac{1}{2(2i)^r\prod_{j=1}^{r}sin(\theta_j)}\sum_{\epsilon_1,..,\epsilon_r\in\{\pm 1 \}}   \prod_{j=1}^r\epsilon_j \sum_{m\in {\mathbb Z}\setminus\{ 0\}}\frac{1}{m^{s+r}}\exp(im\sum_{j=1}^r\epsilon_j\theta_j)(\text{because}\; s\in 2{\mathbb N} ) \\
&=&\frac{1}{2(2i)^r\prod_{j=1}^{r}sin(\theta_j)}\sum_{\epsilon_1,..,\epsilon_r\in\{\pm 1 \}}   (\prod_{j=1}^r\epsilon_j) b_{s+r}(\frac{1}{2\pi}\sum_{j=1}^r\epsilon_j\theta_j)
\end{eqnarray*}
 where $b_k$ is the function defined by $b_k(x)=\sum_ {m\in {\mathbb Z}\setminus\{ 0\}}\frac{e^{2i\pi m x}}{m^k},$ related to the $k$-th Bernouilli polynomial $B_k$ by 
$-\frac{k !}{(2i\pi)^k} b_k(x)=B_k(x-\lfloor x\rfloor)$ with $k\geq 2.$
\medskip

In the case where $G$ is a compact p-adic Lie group, much less are known on these  zeta functions with insertions.
When $p=1$, it is shown (theorem 1.2) in \cite{JZ} that for Fab uniform pro-p group, for fixed $\gamma\in G,$ the function $s\mapsto \zeta^{(1)}(\gamma,s)$ is a rational function in $p^s.$
A question that is pending is the dependence on $\gamma_1,...,\gamma_r.$ There is no result in the litterature  on this question. 
We will answer this question for some selected examples: in the present article  for the case  $GL(2,{\mathbb F}_q), PGL(2,{\mathbb F}_q)$ (section  2) and in \cite{BaRo}  for the case $GL(2,{\mathbb Z}_p), PGL(2,{\mathbb Z}_p)$.

\subsection{Generalized zeta function and non orientable surfaces: Schur-Frobenius formula}

One generalisation  that can be considered  is the case where $\Sigma$ is a compact non orientable surface.  
In this case, $\Sigma$ is homeomorphic to the connected sum of g copies of ${\mathbb P}^2({\mathbb R})$ with $g\geq 1$ and the Euler characteristic of $\Sigma$ is
 $\chi(\Sigma)=2-g.$

We fix a point $x$ on $\Sigma,$
and consider the  set of  group morphisms $X_G(x,\Sigma)=Hom(\pi_1(x,\Sigma),G).$
$X_G(x,\Sigma)/Ad G$ is in bijection with the set of $G$-principal bundle over $\Sigma.$

 A generalisation of 
 Mednykh formula gives \cite{FS}\cite{Sn}:
\begin{equation}
\frac{\vert X_G(x,\Sigma)\vert }{\vert  G\vert}=\vert G\vert^{-\chi(\Sigma)}\sum_{\pi\in \widehat{G},\nu_2(\pi)\not=0}(\nu_2(\pi)
 dim \pi )^{\chi(\Sigma)},
\end{equation}
where $\nu_2(\pi)\in\{+1,0,-1\} $ 
is the Frobenius-Schur indicator of $\pi$ and the sum is restricted to the value of $\pi$ such that $\nu_2(\pi)$ is different of zero. 
The Frobenius indicator is 
 defined  by the formula
\begin{equation}
\label{Frobenius}\nu_2(\pi)=\int_{G}\chi_{\pi}(g^2)dg.
\end{equation}
In fact Frobenius and Schur already  obtained this formula in 1906 while studying the number of solutions  $(R_1,...,R_g)$ in $G^g$ of the equation $R_1^2\cdots R_g^2=e.$

Proof:
This generalisation is easily obtained using the fact that $\pi_1(\Sigma,x)$ admits a presentation by generators $r_1,\cdots, r_g$ with relation $r_1^2\cdots r_g^2=e.$
We have  $X_G(x,\Sigma)=\{(R_1,\cdots, R_g)\in G^g, R_1^2\cdots R_g^2=e\}.$
If we introduce the function  $\theta_{{\mathbb P}^2({\mathbb R})}: G\rightarrow {\mathbb N}, \theta_{{\mathbb P}^2({\mathbb R})}(g)=\vert \{h\in G, h^2=g\}\vert,$ 
it admits the following Fourier decomposition
 \begin{equation}
\label{thetaP2} \theta_{{\mathbb P}^2({\mathbb R})}=\sum_{\pi \in \widehat{G}}\nu_2(\pi) \chi_{\pi}.
 \end{equation}
From 
\begin{eqnarray*}
\vert X_G(x,\Sigma)\vert &=&\vert G\vert^{g} \int_{G^{2g}}   \theta_{{\mathbb P}^2({\mathbb R})}  (k_1)\cdots \theta_{{\mathbb P}^2({\mathbb R})}(k_g)\delta(k_1..k_g)dk_1\cdots dk_g\\
&=&\vert G\vert^{g-1} (\underbrace{\theta_{{\mathbb P}^2({\mathbb R})} *\cdots *\theta_{{\mathbb P}^2({\mathbb R})}}_{g \,\text{times}})(e)
\end{eqnarray*}
we obtain the generalization of Mednikh formula. This is the Schur-Frobenius proof.
$\square$

Therefore we can  define for any rigid compact topological group $G$, the following functions
\begin{equation}
{\zeta}_{G,\epsilon}(s)=\sum_{\pi\in \widehat{G}, \nu_2(\pi)=\epsilon }\frac{1}{(dim \pi)^s}, s\in {\mathbb C},\epsilon\in\{+1,-1,0\}.
\end{equation}
We have ${\zeta}_{G}={\zeta}_{G,+1}+{\zeta}_{G,0}+{\zeta}_{G,-1}$ and $\frac{\vert X_G(x,\Sigma)\vert}{\vert G\vert} =\vert G\vert^{-\chi(\Sigma)}(
{\zeta}_{G,+1}(-\chi(\Sigma))+(-1)^{\chi(\Sigma)}{\zeta}_{G,-1}(-\chi(\Sigma))).$

The value of $ {\zeta}_{G,+1}(n)+(-1)^{n}{\zeta}_{G,-1}(n),$ for $n$  positive integer,   has  been considered in  \cite{Wi}\cite{HoJe} 
 where $G$ is a connected simply connected semisimple compact Lie group, in connection with the volume of the moduli space of $G$ flat connections on $\Sigma,$ where a  relation such as (\ref{volumeofflatconnections})  holds  true.

Mednikh formula has also  a generalisation to the case of compact non orientable surface with boundaries. 

Let $G$ be a finite group, the $\zeta_G$  function can be generalized to the case where one consider connected compact non-orientable  topological surfaces $\Sigma$ with boundaries.
Let $\Sigma$ be a connected compact non-orientable topological surface $\Sigma$ of Euler characteristic $\chi(\Sigma)=2-g$, denote  $\Sigma'= \Sigma \setminus (\cup_{j=1}^r D_j)$ 
where $D_j$ are  open disks included in $\Sigma$  whose closures do not intersect pairwise. We  have $\partial \Sigma'= C_1\cup \cdots \cup C_r,$ where $C_1,..., C_r$ are homeomorphic to circles.
Let $x$ be an interior point of $\Sigma'$, and fix closed curves $\ell_1,...,\ell_r$ originating from $x$ such that the homotopy class of $\ell_j$ is the same as 
$C_j$.
The homotopy group  $\pi_1(\Sigma',x)$ is isomorphic to the group generated by $r_i, i=1,\cdots  g, m_j,  j=1,\cdots r,$ with relation
 $ \prod_{i=1}^g r_i^2\prod_{j=1}^r m_j =e.$

We fix $\gamma_1,\cdots,\gamma_r\in G,$
we can still define the finite set   ${X}_G(x,\Sigma,\gamma_1,...,\gamma_r)$  of principal $G$-bundles over $\Sigma'$  such that the holonomy  along $\ell_j$  belongs to the conjugacy 
class of $\gamma_j,  j=1\cdots r $.

\begin{proposition}
\begin{equation}
\label{}\frac{\vert {X}_G(x,\Sigma,\gamma_1,...,\gamma_r)\vert}{\vert G\vert}=\vert G\vert^{-\chi(\Sigma)}\prod_{j=1}^r\vert O(\gamma_j)\vert
\sum_{\pi\in \widehat{G},\nu_2(\pi)\not=0 }\frac{1}{(\nu_2(\pi)dim(\pi))^{-\chi(\Sigma)}}
\frac{\prod_{j=1}^r \chi_{\pi}(\gamma_j)}{(dim \pi)^{r}}.
\end{equation}
\end{proposition}
Proof:
Direct computation, using
 \begin{eqnarray*}
&&\vert {X}_G(x,\Sigma,\gamma_1,...,\gamma_r)\vert=\\
&&\vert G \vert^{g+r}\int_{G^{g+r}}\delta( \prod_{i=1}^g R_i^2\prod_{j=1}^r M_j)\prod_{j=1}^r\delta_{O(\gamma_j)} (M_j) \prod_i dR_i \prod_j dM_j.
\end{eqnarray*}
$\square.$

We are therefore led to define,  for a  rigid compact  group $G$ and $\epsilon\in \{+1,-1,0 \}$ a  generalized zeta functions $\zeta_{G,\epsilon}^{(r)}:G^r\times
 {\mathbb C}\rightarrow {\mathbb C}$ defined by 
\begin{equation}\zeta_{G,\epsilon}^{(r)}(\gamma_1, ...,\gamma_r) (s)=
\sum_{\pi\in \widehat{G}, \nu_2(\pi)=\epsilon}\frac{\prod_{j=1}^r \chi_{\pi}(\gamma_j)}{(dim \pi)^{s+r}}, s \in {\mathbb C}, \gamma_1,\cdots, \gamma_r\in G.
\end{equation}

To our knowledge these functions do not appear in the litterature and their properties are completely unknown.

\subsection{Generalization of Mednikh formulas and the $\zeta$ function representation of the  double $D(G)$}
Let $G$ be a finite group, $\Sigma$ a  compact  orientable topological surface of genus $g$,  we can express the number $\vert Hom(\pi_1(\Sigma),G)/Ad G\vert $ in term of evaluation at
 $2g-2$ of the 
$\zeta$ representation function of the quantum double of $G.$
This formula appears in \cite{BB} and has been proved using topological field theory. We give here a simple proof of it, answering a question asked in \cite{BB},  which admit straightforward generalisation to the non orientable case and to the case with boundaries that we establish.

Choose a section $\iota :G/ Ad G\rightarrow  G$, and denote $C_\iota(O)=C(\iota(O))$ where $C(k), k\in G$ denotes the centralizer in $G$ of the element $k.$

The double $D(G)$ of a finite group is a quasitriangular  Hopf ${\mathbb C}$-algebra which is the quantum double of the algebra of complex functions on $G.$
 This algebra was introduced and its representations were classified
in \cite{DPR}. The aim of this work was to  understand the work  \cite{DVVV}  from the point of view of braided tensor categories.
It is shown that for each couple $(O,\rho)$ where $O$ is a conjugacy class of $G$ and $\rho$ is an irreducible finite dimensional representation of $C_\iota(O)$ there 
exists an irreducible  finite dimensional representation of $D(G)$ denoted $\Pi{(O,\rho)}$ of dimension $\vert O\vert dim\rho.$
We denote $\widehat{D(G)}$ the set of isomorphism classes of   irreducible representations of $D(G)$.
The set $\{\Pi{(O,\rho)}\}$  are representatives  of $\widehat{D(G)}$  and this is independent of the choice of $\iota.$

We can therefore define the function $\zeta_{D(G)}$
 \begin{eqnarray}
\zeta_{D(G)}(s)&=&\sum_{\Pi\in \widehat{D(G)}}\frac{1}{(dim \Pi)^s}, s\in {\mathbb C}\\
&=&\sum_{O\in G/ Ad G, \rho\in \widehat{C_\iota(O)}}
\frac{1}{(\vert O\vert dim\rho)^s }.
\end{eqnarray}

\begin{proposition}The generalization of Mednikh formula holds:
\begin{eqnarray*}
\vert Hom(\pi_1(\Sigma), G))/Ad G \vert&=&\sum_{k\in G}\frac{\vert  C(k) \vert}{\vert G\vert}\vert C(k)\vert ^{2g-2}\zeta_{C(k)}(2g-2)\\
&=&\sum_{O\in G/Ad G}\vert C_\iota(O)\vert ^{2g-2}\zeta_{C_\iota(O)}(2g-2)\\
&=&\vert G\vert^{2g-2}\zeta_{D(G)}(2g-2).\end{eqnarray*}
\end{proposition}

Proof: $G$ acts on $X_G(x,\Sigma)$ by adjoint action,  applying Burnside counting lemma,  we obtain 
$$
\vert X_G(x,\Sigma)/Ad G\vert=\frac{1}{\vert G\vert}\sum_{k\in G}\vert X_G(x,\Sigma)^k\vert 
$$
where $ X_G(x,\Sigma)^k$ is the subset of elements of $X_G(x,\Sigma)$ fixed by $k.$
From $Hom(.,G)^{k}=Hom(.,C(k)), $ we obtain 
\begin{eqnarray*}
\vert X_G(x,\Sigma)/Ad G\vert&=&\frac{1}{\vert G\vert}\sum_{k\in G}\vert X_{C(k)}(x,\Sigma)\vert \\
&=&\sum_{k\in G}\frac{\vert  C(k) \vert}{\vert G\vert}\vert C(k)\vert ^{2g-2}\zeta_{C(k)}(2g-2).
\end{eqnarray*}

The second equality of the proposition is a consequence of:
\begin{eqnarray*}
&&\sum_{k\in G}\frac{\vert  C(k) \vert}{\vert G\vert}\vert C(k)\vert ^{2g-2}\zeta_{C(k)}(2g-2)=\\
&=&\sum_{k\in G}\frac{1}{\vert O(k)\vert }\sum_{\rho\in \widehat{C(k)}} \vert C(k)\vert^{2g-2}\frac{1}{(dim \rho)^{2g-2}}\\
&=&\sum_{O\in G/Ad G}\sum_{\rho\in \widehat{C_\iota(O)}} \vert C_\iota(O)\vert^{2g-2}\frac{1}{(dim \rho)^{2g-2}}\\
&=& \sum_{O\in G/Ad G}\sum_{\rho\in \widehat{C_\iota(O)}} \vert G\vert^{2g-2}\frac{1}{(\vert O\vert dim \rho)^{2g-2}}\\
&=&\vert G\vert^{2g-2}\zeta_{D(G)}(2g-2).
\end{eqnarray*}
$\square$

These formulas can be generalized to the case of surface with boundaries and non orientable surfaces.

We begin with a lemma.
\begin{lemma}
Let $H$ be a subgroup of $G$ and $\gamma\in G$, the restriction to $H$ of the characteristic function of the conjugacy class $O(\gamma)$ of $G$ admits the following Fourier transform:
\begin{equation*}
\delta_{O(\gamma)}\vert_{H}=\frac{1}{\vert C(\gamma)\vert}\sum_{\rho\in \widehat{H}}\chi_{\rho}Tr(Ind_{H}^G(\rho)(\gamma^{-1})).
\end{equation*}
\end{lemma}
Proof:
We use the classical notation of group theory: if $\phi$ is a function on $H$ we extend it on $G$ as $\phi^0$ with the definition $\phi^0(x)=\phi(x)$ if $x\in H$ and zero otherwise.
The Fourier expansion of the class function $\delta_{O(\gamma)}\vert_{H}$ is $\sum_{\rho\in \widehat{H}}a_\rho\chi_\rho,$
therefore:
\begin{eqnarray*}
\overline{a}_\rho&=&\frac{1}{\vert H\vert}\sum_{h\in H}\delta_{O(\gamma)}(h)\chi_\rho(h)\\
&=&\frac{1}{\vert H\vert}\sum_{h\in H\cap O(\gamma)}\chi_\rho(h)\\
&=&\frac{1}{\vert H\vert}\sum_{ h\in O(\gamma)}\chi_\rho^{0}(h)\\
&=&\frac{1}{\vert H\vert\vert C(\gamma)\vert}\sum_{x\in G}\chi_\rho^{0}(x\gamma x^{-1})\\
&=&\frac{1}{\vert C(\gamma)\vert} Tr(Ind_{H}^G(\rho)(\gamma)).
\end{eqnarray*}
$\square$

\begin{proposition}The generalization of Mednikh formula holds:
\begin{eqnarray*}
&&\vert X_G(\Sigma,\gamma_1,..,\gamma_r)/Ad G \vert=\\
&=&\frac{\prod_{j=1}^r\vert O(\gamma_j)\vert}{\vert G\vert^{r+1} }
\sum_{k\in G}\vert C(k)\vert ^{2g+r-1}\sum_{\rho\in\widehat{C(k)}}\frac{\prod_{j=1}^r Tr(Ind_{C(k)}^G(\rho)(\gamma_j))}{(dim \rho)^{2g-2+r}}
\end{eqnarray*}
\end{proposition}
Proof: Straightforward
$\square.$

In the case where $\Sigma$ is a compact non orientable surface of Euler characteristic $2-g$ we have:

\begin{proposition}
\begin{equation*}
\vert Hom(\pi(\Sigma), G))/Ad G \vert=\sum_{O\in G/Ad G}\vert C_\iota(O)\vert^{-\chi(\Sigma)}\sum_{\rho\in \widehat{ C_\iota(O)},\nu_2(\rho)\not=0}(\nu_2(\rho)
 dim \rho )^{\chi(\Sigma)},
\end{equation*}
\end{proposition}
Proof: Straightforward.
$\square$

In the case where $\Sigma$ is a compact non orientable surface of Euler characteristic $2-g$, with boundaries we have:

\begin{proposition}
\begin{eqnarray*}
&&\vert X_G(\Sigma,\gamma_1,..,\gamma_r)/Ad G \vert=\\
&=&\frac{\prod_{j=1}^r\vert O(\gamma_j)\vert}{\vert G\vert^{r+1} }
\sum_{k\in G}\vert C(k)\vert^{g+r-1}\sum_{\rho\in\widehat{C(k)},\nu_2(\rho)\not=0}\frac{\prod_{j=1}^r Tr(Ind_{C(k)}^G(\rho)(\gamma_j))}{(\nu_2(\rho)dim \rho)^{-\chi(\Sigma)}(\dim \rho)^r}
\end{eqnarray*}
\end{proposition}
Proof:Straightforward.
$\square.$

We will now explicitely compute the functions defined in this section in two simple cases $GL(2,{\mathbb F}_q)$ and $PGL(2,{\mathbb F}_q).$

\section{Explicit computation for the groups $GL(2,{\mathbb F}_q)$, $PGL(2,{\mathbb F}_q)$}

\subsection
{Computation of generalized $\zeta$ functions in the $GL(2,{\mathbb F}_q)$ case}
We have collected in the appendix the description of conjugacy classes, irreducible representations of these groups and values of the Frobenius-Schur indicator.
We use the notation ${\mathbb F}={\mathbb F}_q$ and ${\mathbb E}={\mathbb F}_{q^2}.$

From the  classification of irreducible representations of $GL(2,{\mathbb F}_q)$, the zeta function is straightforwardly computed and reads:
\begin{equation}
\zeta_{GL(2,{\mathbb F}_q)} (s)=q-1+\frac{(q-1)(q-2)}{2(q+1)^{s}}+\frac{q-1}{q^{s}}+\frac{(q-1)q}{2(q-1)^{s}}.
\end{equation}

Let $\gamma_1,\cdots, \gamma_r$ elements of $G=GL(2,{\mathbb F}_q),$ we consider $\zeta_G^{(r)}(\gamma_1,\cdots,\gamma_r)(s).$ 
Note that the group of characters of $G$ acts by tensorisation on the set of irreducible representations, therefore we have, for any $\mu\in \widehat{{\mathbb F}_q^\times}$
\begin{equation}
\zeta_G^{(r)}(\gamma_1,\cdots,\gamma_r)(s)=\chi_{\mu}(\gamma_1\cdots\gamma_r)\zeta_G^{(r)}(\gamma_1,\cdots,\gamma_r)(s).
\end{equation}
This implies that  the function $\zeta_G^{(r)}(\gamma_1,\cdots,\gamma_r)(s)$ is null if $det(\gamma_1\cdots\gamma_r)\not=1.$

We will  compute the generalized  zeta functions $\zeta^{(1)} $ for the 4 types of conjugacy classes.

These can easily be done using the  identities valid for any $x\in {\mathbb F}^{\times}$
\begin{eqnarray*}
&&\sum_{\mu\in  \widehat{{\mathbb F}^{\times}} } \mu(x)=(q-1)\delta(x),\\
&&\sum_{\substack{\{\mu_1,\mu_2\}\subset  \widehat{{\mathbb F}^{\times}}\\\mu_1\not=\mu_2}}\mu_1(x)\mu_2(y)= 
\frac{1}{2}(\sum_{ \mu\in \widehat{{\mathbb F}^{\times}}}\mu(x))^2-\frac{1}{2}\sum_{ \mu\in \widehat{{\mathbb F}^{\times}}}\mu(x)^2=\frac{1}{2}(q-1)^2\delta({x})-\frac{1}{2}(q-1)
\delta(x^2),\\
&&\sum_{\nu \in  (\widehat{{\mathbb E}^\times}\setminus  \widehat{{\mathbb F}^\times})/{\sim}}\nu(x)=
\frac{1}{2}\sum_{\nu \in \widehat{\mathbb E}^\times}\nu(x) -\frac{1}{2}\sum_{\mu \in \widehat{\mathbb F}^\times}\mu(N(x)) =\frac{1}{2}(q^2-1)\delta(x)- \frac{1}{2}(q-1)\delta({x^2}).
\end{eqnarray*}

We therefore have 
$$\zeta^{(1)}(c_1(x))(s)=\sum_{\mu\in \widehat{{\mathbb F}^{\times} }}\frac{\mu(x)^2}{1^{s+1}}+
\sum_{\substack{\{\mu_1,\mu_2\}\subset  \widehat{{\mathbb F}^{\times}}\\\mu_1\not=\mu_2}}\frac{(q+1)\mu_1(x)\mu_2(x)}{(q+1)^{s+1}}+
\sum_{\mu\in \widehat{ {\mathbb F}^{\times} }}\frac{q\mu(x)^2}{q^{s+1}}+\sum_{\nu\in (\widehat{{\mathbb E}^\times}\setminus  \widehat{{\mathbb F}^\times})/{\sim}}\frac{(q-1)\nu(x)}{(q-1)^{s+1}}$$
$$=(q-1)\delta(x^2)+\frac{1}{2}(q+1)^{-s}((q-1)^2\delta(x)-(q-1)\delta(x^2))+$$ $$q^{-s}(q-1)\delta(x^2)+
\frac{1}{2}(q-1)^{-s}((q^2-1) \delta(x)-(q-1)\delta(x^2)).$$

$$\zeta^{(1)}(c_2(x))(s)=\sum_{\mu\in \widehat{{\mathbb F}^{\times} }}\frac{\mu(x)^2}{1^{s+1}}+
\sum_{\substack{\{\mu_1,\mu_2\}\subset\widehat{ {\mathbb F}^{\times}} \\\mu_1\not=\mu_2}}\frac{\mu_1(x)\mu_2(x)}{(q+1)^{s+1}}
+\sum_{\nu \in (\widehat{{\mathbb E}^\times}\setminus  \widehat{{\mathbb F}^\times})/{\sim}}\frac{-\nu(x)}{(q-1)^{s+1}}$$
$$=(q-1)\delta(x^2)+\frac{1}{2}(q+1)^{-s-1}((q-1)^2\delta(x)-(q-1)\delta(x^2))-\frac{1}{2}(q-1)^{-s-1}((q^2-1)\delta(x)-(q-1)\delta(x^2)).$$

$$\zeta^{(1)}(c_3(x,y))(s)=\sum_{\mu\in \widehat{{\mathbb F}^{\times}} }\frac{\mu(x)\mu(y)}{1^{s+1}}+
\sum_{\substack{\{\mu_1,\mu_2\}\subset\widehat{ {\mathbb F}^{\times} }\\\mu_1\not=\mu_2}}\frac{\mu_1(x)\mu_2(y)+\mu_1(y)\mu_2(x)}{(q+1)^{s+1}}
+\sum_{\mu\in\widehat{ {\mathbb F}^{\times} }}\frac{\mu(xy)}{q^{s+1}}$$
$$=(q-1)\delta(xy)+(q+1)^{-s-1}((q-1)^2\delta(x)\delta(y)-(q-1)\delta(xy))+q^{-s-1}(q-1)\delta(xy).$$

$$\zeta^{(1)}(c_4(\lambda))(s)=\sum_{\mu\in \widehat{{\mathbb F}^{\times} }}\frac{\mu(\lambda\overline{\lambda})}{1^{s+1}}-
\sum_{\mu \in \widehat{{\mathbb F}^{\times}}}\frac{\mu(\lambda\overline{\lambda})}{q^{s+1}}-\sum_{\nu \in (\widehat{{\mathbb E}^\times}\setminus  \widehat{{\mathbb F}^\times})/{\sim}}\frac{\nu(\lambda)+\nu(\overline{\lambda})}{(q-1)^{s+1}}$$
$$=(q-1)\delta({\lambda\overline{\lambda}})-(q-1)q^{-s-1}\delta({\lambda\overline{\lambda}})-(q^2-1)(q-1)^{-s-1}\delta({\lambda})+(q-1) (q-1)^{-s-1}\delta({\lambda\overline{\lambda}}).$$

  $\zeta^{(1)}(\gamma, s)$ vanishes when the determinant of $\gamma$ is not equal to $1,$ 
this is in sharp contrast with the $PGL(2,{\mathbb F}_q)$ case  that we will consider in the next section.
\medskip

We can compute explicitely the zeta function of the double $D(G).$

We first describe the centralizer of elements of $G.$

\begin{itemize}
\item 
$C(c_1(x))=GL(2,{\mathbb F}).$
\item
$C(c_2(x))=\{\begin{pmatrix}a&b\\0&a
\end{pmatrix}, a\in {\mathbb F}^\times, b\in {\mathbb F}\}\simeq {\mathbb F}^\times\times {\mathbb F}$ abelian group of order $q(q-1).$
\item 
$C(c_3(x,y))=\{\begin{pmatrix}a&0\\0&b
\end{pmatrix}, (a,b)\in ({\mathbb F}^\times)^2\}\simeq ({\mathbb F}^\times)^2$ abelian group of order $(q-1)^2.$
\item 
$C(c_4(\lambda))=\{\begin{pmatrix}a&-c\lambda\overline{\lambda}\\c&a+c(\lambda+\overline{\lambda})
\end{pmatrix}, (a,b)\in  {\mathbb F}^2\setminus \{(0,0)\}\}$ abelian group of order $q^2-1.$
\end{itemize}

Therefore we have 
\begin{eqnarray*}
&&\zeta_{D(GL(2,{\mathbb F}))}(s)=\sum_{x\in {\mathbb F}^\times}\zeta_{C(c_1(x))}(s)+
\sum_{x\in {\mathbb F}^\times}\frac{1}{\vert C(c_2(x))\vert^s}\vert  {\mathbb F}^\times\times {\mathbb F}\vert+\\
&&+\sum_{\substack{\{ x,y\} \subset {\mathbb F}^\times\\ x \not=y}} \frac{1}{\vert C(c_3(x,y))\vert^s }\vert  {\mathbb F}^\times\times  {\mathbb F}^\times\vert +
\sum_{\substack{\lambda \in  {\mathbb E}\setminus  {\mathbb F}\\ \lambda\sim \overline{\lambda}} }\frac{1}{\vert C(c_4(\lambda))\vert^s }\vert  {\mathbb F}^2\setminus \{(0,0)\}\vert.
\end{eqnarray*}

\begin{eqnarray*}
&&\zeta_{D(GL(2,{\mathbb F}))}(s)=(q-1)^2+\frac{(q-1)^2(q-2)}{2(q+1)^{s}}+\frac{(q-1)^2}{q^{s}}+\frac{(q-1)^2q}{2(q-1)^{s}}+\\
&&+\frac{q(q-1)^2}{(q(q-1))^{s}}+\frac{(q-1)^3(q-2)}{2(q-1)^{2s}}+\frac{q(q-1)^2(q+1)}{2(q^2-1)^{s}}.
\end{eqnarray*}

We now give the form of 
 $\zeta_{GL(2,{\mathbb F}_q),\epsilon}$ functions.

There is no Frobenius Schur indicator having $-1$ value, therefore  $\zeta_{GL(2,{\mathbb F}),-1}=0.$

When $q$ is odd  we have  $$\zeta_{GL(2,{\mathbb F}),1}(s)=2+\frac{q-1}{2(q+1)^s}+
\frac{2}{q^s}+\frac{q-1}{2(q-1)^s},$$
when $q$ is even we have $$\zeta_{GL(2,{\mathbb F}),1}(s)=1+\frac{q-2}{2(q+1)^s}+
\frac{1}{q^s}+\frac{q}{2(q-1)^s}.$$

The expression of $\zeta_{GL(2,{\mathbb F}),0}$ is obtained from $\zeta_{GL(2,{\mathbb F}),0}=\zeta_{GL(2,{\mathbb F})}-\zeta_{GL(2,{\mathbb F}),1}.$

\subsection
{Computation of generalized $\zeta$ functions in the $PGL(2,{\mathbb F}_q)$ case}

From the structure classification of irreducible representations of $PGL(2,{\mathbb F}_q)$, the zeta function is straightforwardly computed and reads:
\begin{itemize}
\item In odd characteristic one has:
$\zeta_{PGL(2,{\mathbb F}_q)} (s)=2+\frac{(q-3)}{2(q+1)^s}+\frac{2}{q^s}+\frac{(q-1)}{2(q-1)^s}$
\item in characteristic $2$ one has:
$\zeta_{PGL(2,{\mathbb F}_q)} (s)=1+\frac{(q-2)}{2(q+1)^s}+\frac{1}{q^s}+\frac{q}{2(q-1)^s}.$
\end{itemize}

From the property of the Schur-Frobenius indicator of irreducible representations of
 $PGL(2,{\mathbb F}),$ we have:
$$\zeta_{PGL(2,{\mathbb F}),-1}=\zeta_{PGL(2,{\mathbb F}),0}=0.$$
\medskip 

We now give explicit expressions of zeta functions with insertions.

We denote by ${\mathcal M}$ the set of equivalence classes for  the relation $\mu\sim \mu^{-1}$ on the set of  characters     $\mu\in  
\widehat{{\mathbb F}^{\times}} $ such that     $\mu\not=1$ in the case $q$ even and $\mu\notin\{1,\epsilon\}$  in the case $q$ odd. This set is in bijection with the set of inequivalent representations of 
$PGL(2,{\mathbb F})$ of principal  type.

We denote by ${\mathcal N}$ the set of equivalence classes for  the relation $\nu\sim \nu^{-1}$ on the set of primitive characters     $\nu\in \widehat{{\mathbb E}^{\times}} \setminus 
\widehat{{\mathbb F}^{\times} }$ such that $\nu\vert_{{\mathbb F}^{\times} }=id.$  This set is in bijection with the set of inequivalent representations of 
$PGL(2,{\mathbb F})$ of cuspidal type.

$\bullet$ We first  assume that  $q$ is odd.

Computation of the  $\zeta$ function with one insertion.

$$\zeta^{(1)}(\begin{pmatrix}1&1\\0&1
\end{pmatrix},s)=2+\frac{(q-3)}{2(q+1)^{s+1}}-\frac{(q-1)}{2q^{s+1}}.$$
Let $x\in {\mathbb F}^\times,x\not=1,$
$$\zeta^{(1)}(\begin{pmatrix}x&0\\0&1
\end{pmatrix},s)=1+\epsilon(x)+\sum_{\mu\in  {\mathcal M}}\frac{\mu(x)+\mu(x)^{-1}}{(q+1)^{s+1}}+\frac{1+\epsilon(x)}{q^{s+1}}$$
$$=1+\epsilon(x)+\frac{1}{(q+1)^{s+1}}((q-1)\delta(x)-1-\epsilon(x))+\frac{1+\epsilon(x)}{q^{s+1}}$$
$$=(1+\frac{1}{q^{s+1}}-\frac{1}{(q+1)^{s+1}})(1+\epsilon(x)).$$

Let $\epsilon_{\mathbb E}:{\mathbb E}^{\times}\rightarrow \{-1,1\}$ be the morphism of multiplicative group defined by  $\epsilon_{\mathbb E}(\lambda)=-1$ if and only if  $\lambda$ is not a square in ${\mathbb E}^{\times}.$

Let $\lambda\in {\mathbb E}\setminus {\mathbb F}$,
$$\zeta^{(1)}(\begin{pmatrix}0&-\lambda\overline{\lambda}\\1&\lambda+\overline{\lambda}
\end{pmatrix},s)=$$
$$=1+\epsilon(\lambda\overline{\lambda})-\sum_{\nu\in {\mathcal N}}\frac{\nu(\lambda)+\nu(\overline{\lambda})}{(q-1)^{s+1}}-\frac{1+\epsilon(\lambda\overline{\lambda})}{q^{s+1}}$$
$$=(1+\frac{1}{(q-1)^{s+1}}-\frac{1}{q^{s+1}})+\epsilon(\lambda\overline{\lambda})(1-\frac{1}{q^{s+1}})+\frac{1}{(q-1)^{s+1}}\epsilon_{\mathbb E}(\lambda).$$
where we have used the identity $\sum_{\nu\in {\mathcal N}}(\nu(\lambda)+\nu(\overline{\lambda}))=-1-\epsilon_{\mathbb E}(\lambda),$ for $\lambda\in {\mathbb E}\setminus {\mathbb F}.$

Computation of the generalized  function  with r-insertions in the case where all the conjugacy classes are diagonal.

Let $\gamma_i=\begin{pmatrix}x_i&0\\0&1
\end{pmatrix}, i=1,..,r$ with $x_i\in {\mathbb F}^{\times}\setminus\{1\},$ we have 
$$\zeta^{(r)}(\gamma_1,\cdots,\gamma_r,s)=1+\epsilon(x_1\cdots x_r)+\sum_{\mu\in  {\mathcal M} }\frac{(\mu(x_1)+\mu(x_1)^{-1})\cdots (\mu(x_r)+\mu(x_r)^{-1})}{(q+1)^{s+r}}+$$
 $$+\frac{1+\epsilon(x_1\cdots x_r)}{q^{s+r}}$$
$$=(1+\epsilon(x_1\cdots x_r))(1+\frac{1}{q^{s+r}})+
\frac{1}{2(q+1)^{s+r}}\sum_{\substack{\eta_i\in\{\pm 1\}\\i=1,..,r}}((q-1)\delta({x_1^{\eta_1}\cdots x_r^{\eta_r}})-1-\epsilon(x_1^{\eta_1}\cdots x_r^{\eta_r})).$$

Computation of the generalized  function  with r-insertions in the case where all the conjugacy classes are elliptic.

Let $\gamma_i=\begin{pmatrix}0&-\lambda_i\overline{\lambda_i}\\1&\lambda_i+\overline{\lambda_i}
\end{pmatrix}, i=1,..,r$ with $\lambda_i\in {\mathbb E}^{\times}\setminus{\mathbb F},$ we have 
$$\zeta^{(r)}(\gamma_1,\cdots,\gamma_r,s)=1+\epsilon(\lambda_1\overline{\lambda}_1\cdots\lambda_r\overline{\lambda}_r )+(-1)^r\sum_{\nu\in {\mathcal N}   }\frac{(\nu(\lambda_1)+\nu(\overline{\lambda_1}))\cdots (\nu({\lambda_r})+\nu(\overline{\lambda}_r))}{(q-1)^{s+r}}+
 (-1)^r\frac{1}{q^{s+r}}$$
$$=(1+(-1)^r\frac{1}{q^{s+r}})+
\frac{(-1)^{r}}{2(q-1)^{s+r}}\sum_{\substack{\eta_i\in\{id, \sigma\}\\i=1,..,r}}
((q+1)\phi_{\mathbb F}(
\eta_1(\lambda_1)..\eta_r(\lambda_r))-1-\epsilon_{\mathbb E}(\eta_1(\lambda_1)..\eta_r(\lambda_r))),
$$

where $\phi_{\mathbb F}:{\mathbb E}\rightarrow\{0,1\}$ is the characteristic function of ${\mathbb F}$ which satisfies 
$\sum_{\nu\in {\mathcal N}}(\nu(\lambda)+\nu(\overline{\lambda}))=(q+1)\phi_{\mathbb F}(\lambda)-1-\epsilon_{\mathbb E}(\lambda)$ for $\lambda\in {\mathbb E}^{\times}.$
\medskip

Computation of the generalized  function  with r insertions in the mixed case where m are diagonal and n are elliptic, with $m\geq 1$ and $n\geq 1.$

We denote $\gamma_i=\begin{pmatrix}x_i&0\\0&1
\end{pmatrix}, i=1,..,m, $ with $x_i\in {\mathbb F}^{\times}\setminus\{1\},$ $\gamma_{i+m}=\begin{pmatrix}0&-\lambda_i\overline{\lambda_i}\\1&\lambda_i+\overline{\lambda_i}
\end{pmatrix}, i=1,..,n$ with $\lambda_i\in {\mathbb E}^{\times}\setminus{\mathbb F},$ the structure of the character table implies that 
$$\zeta^{(r)}(\gamma_1,\cdots,\gamma_{m+n},s)=
(1+\frac{(-1)^n}{q^{s+r}})(1+\epsilon(x_1..x_m\lambda_1\overline{\lambda_1}...\lambda_n\overline{\lambda_n}) ).$$

$\bullet$ We now assume   $q$ is a power of $2.$

Computation of the  $\zeta$ function with one insertion.

$$\zeta^{(1)}(\begin{pmatrix}1&1\\0&1
\end{pmatrix},s)=1+\frac{(q-2)}{2(q+1)^{s+1}}-\frac{(q)}{2q^{s+1}}.$$
Let $x\in {\mathbb F}^\times,x\not=1,$
$$\zeta^{(1)}(\begin{pmatrix}x&0\\0&1
\end{pmatrix},s)=1+\sum_{\mu\in { \mathcal M} }\frac{\mu(x)+\mu(x)^{-1}}{(q+1)^{s+1}}+
\frac{1}{q^{s+1}}$$
$$=1+\frac{1}{(q+1)^{s+1}}((q-1)\delta(x)-1)+\frac{1}{q^{s+1}}=1-\frac{1}{(q+1)^{s+1}}+\frac{1}{q^{s+1}}. $$

Let $\lambda\in {\mathbb E}\setminus {\mathbb F}$,
$$\zeta^{(1)}(\begin{pmatrix}0&-\lambda\overline{\lambda}\\1&\lambda+\overline{\lambda}
\end{pmatrix},s)=1-\sum_{\nu\in {\mathcal N}}\frac{\nu(\lambda)+\nu(\overline{\lambda})}{(q-1)^{s+1}}-\frac{1}{q^{s+1}}$$
$$=1-\frac{1}{q^{s+1}}+\frac{1}{(q-1)^{s+1}},$$
where we have used the identity $\sum_{\nu\in {\mathcal N}}(\nu(\lambda)+\nu(\overline{\lambda}))=-1$ for $\lambda\in {\mathbb E}\setminus {\mathbb F}.$

Computation of the generalized  function  with r insertions in the case where all the conjugacy classes are diagonal.

Let $\gamma_i=\begin{pmatrix}x_i&0\\0&1
\end{pmatrix}, i=1,..,r$ with $x_i\in {\mathbb F}^{\times}\setminus\{1\},$ we have 
$$\zeta^{(r)}(\gamma_1,\cdots,\gamma_r,s)=1+\sum_{\mu\in{\mathcal M}}\frac{(\mu(x_1)+\mu(x_1)^{-1})\cdots (\mu(x_r)+\mu(x_r)^{-1})}{(q+1)^{s+r}}+\frac{1}{q^{s+r}}$$
$$=(1+\frac{1}{q^{s+r}})+
\frac{1}{2(q+1)^{s+r}}\sum_{\eta_i=\pm 1}((q-1)\delta(x_1^{\eta_1}\cdots x_p^{\eta_r})-1).$$

Computation of the generalized  function  with r insertions in the case where all the conjugacy classes are elliptic.

Let $\gamma_i=\begin{pmatrix}0&-\lambda_i\overline{\lambda_i}\\1&\lambda_i+\overline{\lambda_i}
\end{pmatrix}, i=1,..,r$ with $\lambda_i\in {\mathbb E}^{\times}\setminus{\mathbb F},$ we have 
$$\zeta^{(r)}(g_1,\cdots,g_r,s)=1+(-1)^r\sum_{\nu\in {\mathcal N}   }\frac{(\nu(\lambda_1)+\nu(\overline{\lambda_1}))\cdots (\nu({\lambda_r})+\nu(\overline{\lambda}_r))}{(q-1)^{s+r}}\frac{1}{q^{s+r}}$$
$$=(1-\frac{1}{q^{s+r}})+\frac{(-1)^{r+1}}{2(q-1)^{s+r}}\sum_{\eta_i\in\{id, \sigma\}}
((q+1)\phi_{\mathbb F}(
\eta_1(\lambda_1)..\eta_r(\lambda_r))-1),
$$

where $\phi_{\mathbb F}:{\mathbb E}\rightarrow\{0,1\}$ is the characteristic function of ${\mathbb F}$ and satisfies 
$\sum_{\nu\in {\mathcal N}}(\nu(\lambda)+\nu(\overline{\lambda}))=(q+1)\phi_{\mathbb F}(\lambda)-1$ for $\lambda\in {\mathbb E}^{\times}.$

Computation of the generalized  function  with r insertions in the mixed case where m are diagonal and n are elliptic, with $m\geq 1$ and $n\geq 1.$

We denote $\gamma_i=\begin{pmatrix}x_i&0\\0&1
\end{pmatrix}, i=1,..,m, $ with $x_i\in {\mathbb F}^{\times}\setminus\{1\},$ $\gamma_{i+m}=\begin{pmatrix}0&-\lambda_i\overline{\lambda_i}\\1&\lambda_i+\overline{\lambda_i}
\end{pmatrix}, i=1,..,n$ with $\lambda_i\in {\mathbb E}^{\times}\setminus{\mathbb F},$ the structure of the character table implies that 
$$\zeta^{(r)}(\gamma_1,\cdots,\gamma_r,s)=1+(-1)^n\frac{1}{q^{s+r}}.$$

\section{Conclusion}
This article is giving  a review of Mednikh formulas and its generalisations  and  is pushing  forward the definition of generalized $\zeta$ representation functions of groups.
These generalized functions have not received up to now any attention in the case of p-adic Lie groups, although they are under active study in the case of compact real groups.

The results obtained in this article, although very elementary, are the first necessary step to understand generalized zeta functions with insertions for groups such as 
$GL(2,{\mathcal O})$ and $PGL(2,{\mathcal O})$  which will be developped in \cite{BaRo}.
There are different lines of investigation that can be pursued.

It would be interesting to understand the structure of $\zeta_{G,\epsilon}$ for $G$ p-adic Lie group. 
Can we apply the same techniques developped in \cite{JZ}?

Is it possible to have simple closed expression for  the computation of the zeta function generalisation in the case of $GL(n,{\mathbb F_q})$, $PGL(n,{\mathbb F_q})$, and more generally in the case  $GL(n, {\mathcal O}), PGL(n, {\mathcal O})?$

Can we generalize the notion of quantum double to the case of compact p-adic Lie groups along the line for example of \cite{KM} and define their  generalized $\zeta$ function representations?

\section{Acknowledgments}
 I thank   I.Badulescu for discussions.

\begin{appendix}
\section{Table of characters , Schur Frobenius indicator  and explicit decomposition of tensor products of irreducible representations of the groups 
$GL(2,{\mathbb F}_q)$, $PGL(2,{\mathbb F}_q)$.}
\subsection{Generalities}

We have used  the work of  J.Adams\cite{Ad} and we have completed it.

Let ${\mathbb F}={\mathbb F}_q$  finite field with $q=p^e$ element where $p$ is any prime number and $e\geq 1.$
We have $GL(2,{\mathbb F})=\{A\in Mat(2, {\mathbb F}), det(A)\in {\mathbb F}^{\times}\}$. The center is equal to $Z=Z(GL(2,{\mathbb F}))={\mathbb F}^\times I,$  where $I$ is the unit matrix. We denote $PGL(2,{\mathbb F})=GL(2,{\mathbb F})/Z.$ We have  $SL(2,{\mathbb F})=
\{A\in Mat(2, {\mathbb F}), det(A)=1\}.$ The centre $Z(SL(2,{\mathbb F}))=\{ \lambda I,  \lambda\in {\mathbb F}, \lambda^2=1\},$
and $PSL(2,{\mathbb F})=SL(2,{\mathbb F})/Z( SL(2,{\mathbb F})).$
Depending on the parity of $p$, we have $Z( SL(2,{\mathbb F}))=\{I\}$ if $p=2$ and $Z( SL(2,{\mathbb F})=\{I,-I\}$  if $p$ is odd.

We have $\vert GL(2, {\mathbb F})\vert=q(q-1)^2(q+1), $
$\vert PGL(2, {\mathbb F})\vert=q(q-1)(q+1), $$ \vert SL(2, {\mathbb F})\vert=q(q-1)(q+1)$ and 
$\vert PSL(2, {\mathbb F})\vert=\frac{q(q-1)(q+1)}{gcd(2,q-1)}.$
As a result when $q$ is odd all these groups are different, when $q$ is even but different of $2$ we have 
 $GL(2,{\mathbb F})\not=PSL(2,{\mathbb F})=PGL(2,{\mathbb F})=SL(2,{\mathbb F}).$
When $q=2$ we have $GL(2,{\mathbb F})=PSL(2,{\mathbb F})=PGL(2,{\mathbb F})=SL(2,{\mathbb F})=S_3.$

We denote ${\mathbb E}={\mathbb F}_{q^2},$  extension of ${\mathbb F}$ of dimension $2$ and denote $\sigma$ the generator of $Gal({\mathbb E}, {\mathbb F})={\mathbb Z}/2{\mathbb Z}$. $\sigma$ is the Frobenius element, we have $\sigma(a)=a^q,\forall a\in {\mathbb E},$ and because ${\mathbb E}$ is the unique quadratic extension of ${\mathbb F}$ it will be convenient to denote  $\sigma(a)=\overline{a}.$ As usual we  denote the norm $N$ and the trace $Tr$ being  the maps from ${\mathbb E}$ to ${\mathbb F}$ defined by $N(\lambda)=\lambda \overline{\lambda}=\lambda^{q+1}$ and $Tr(\lambda)=\lambda+\overline{\lambda}=\lambda+\lambda^q.$

We have to consider the two cases $q$ even and $q$ odd.

When $q$ is odd  the map ${\mathbb F}\rightarrow {\mathbb F}, x\mapsto x^2$ is not surjective. Let $\Delta$ a non zero element  of  ${\mathbb F}$ which is not a square. 
Let $\delta\in {\mathbb E}$ a square root of $\Delta.$ We have ${\mathbb E}={\mathbb F}(\delta)$, and $\overline{\delta}=\delta^q=-\delta.$
As a result if $x\in {\mathbb E}$ we can uniquely write $x=a+b\delta$ with $a,b\in {\mathbb F}$ and $Tr(x)=2a$ and $N(x)=a^2-b^2\Delta.$

When $q$ is a power of $2$,  every element of ${\mathbb F}$ is a  square. Let us consider instead  the map ${\mathbb F}\rightarrow {\mathbb F}, x\mapsto x^2+x.$
This is a ${\mathbb F}_2$ linear  map, never surjective and  its image  has $\frac{1}{2}q$ elements. Let $\Omega$ be any element which is not in the image of this map.
We have  $ {\mathbb E}= {\mathbb F}(\omega)$ where $\omega $ is a root of the  irreducible polynomial  $X^2+X+\Omega.$

\subsection{Character table of $GL(2,{\mathbb F})$}
\subsubsection{Conjugacy classes}

We first give a list of  the conjugacy classes, which are classified into $4$ types, the number of conjugacy classes of this type and the number of element of each conjugacy class.
For each of them we choose a particular representative.
Two objects  will be in relation $\backsim$ if and only if their corresponding conjugacy class are equal.

\begin{itemize}
\item (Central) $c_1(x)=\begin{pmatrix}x&0\\0&x
\end{pmatrix}
, x\in {\mathbb F}^{\times}$\\
Number of classes: $q-1$, Number of elements in each conjugacy classes: 1
\item (Unipotent) $c_2(x)=\begin{pmatrix}x&1\\0&x
\end{pmatrix} ,x\in {\mathbb F}^{\times}$\\
Number of classes: $q-1$, Number of elements in each conjugacy classes: $(q-1)(q+1)$
\item (Diagonal) $ c_3(x,y)=\begin{pmatrix}x&0\\0&y
\end{pmatrix}, (x,y)\in ({\mathbb F}^{\times})^2, x\not=y, (x,y)\backsim(y,x)$\\
Number of classes: $\frac{(q-1)(q-2)}{2}$, Number of elements in each conjugacy classes: $q(q+1)$
\item (Elliptic) $c_4(\lambda)=\begin{pmatrix}0&-\lambda\overline{\lambda}\\1&\lambda+\overline{\lambda}
\end{pmatrix},\lambda\in {\mathbb E}\setminus{\mathbb F}, \lambda\backsim \overline{\lambda}.$\\Number of classes: $\frac{q(q-1)}{2}$, Number of elements in each conjugacy classes: $q(q-1)$

\end{itemize}

\subsubsection{Irreducible representations}

We now give the list of finite dimensional complex irreducible representations of $GL(2,{\mathbb F}).$
 The group of  one dimensional complex  representations of ${\mathbb F}^\times$ is the group of multiplicative characters and is denoted $\widehat{{\mathbb F}^\times}.$ 

We use the implicit definition of a relation  $\sim$: two objects in a set will be in relation   $\sim$  if and only if their corresponding representation are equivalent.

There are 4 types of irreducible representations.
\begin{itemize}
\item Let $\mu$ be a multiplicative character of ${\mathbb F}^\times,$ we denote $\chi_\mu$ the one dimensional representation of $G=GL(2,{\mathbb F})$ defined by $\chi_\mu(g)=\mu(det(g)), \forall g\in G.$ These are the one dimensional representations of $G$, i.e the characters of $G.$ 
There are $q-1$ inequivalent representations of this type.

\item Let $T$ be the diagonal torus and $N$ the unipotent group of upper triangular elements,  we denote $B$ the Borel subgroup, we have $B=TN.$
Let $\mu_1,\mu_2$ be multiplicative characters, we  define a one dimensional representations $\mu_1\otimes \mu_2$ of $T$, extend it  trivially to $N$ to a one  dimensional representation still denoted $\mu_1\otimes \mu_2$  of $B$. We denote $I(\mu_1,\mu_2)=Ind_B^G(\mu_1\otimes \mu_2).$ These representations are of dimension $[G:B]=q+1.$ They  are irreducible if and only if $\mu_1\not=\mu_2.$ We have $I(\mu_1,\mu_2)\simeq I(\mu_2,\mu_1),$  therefore $ (\mu_1,\mu_2)\sim(\mu_2,\mu_1).$ These representations are the principal representations. There are  $\frac{1}{2}(q-1)(q-2)$ inequivalent representations 
of this type.

\item The representation $I(\mu,\mu)$ is not irreducible, it contains $\chi_\mu$ as a subrepresentation. The representation $St_\mu=I(\mu,\mu)/\chi_\mu$ is irreducible of dimension $q$. For $\mu=1$, the representation $St_1$ is usually called Steinberg representation and denoted $St$. We have $St_\mu\simeq \chi_\mu\otimes St. $ There are $q-1$ inequivalent representations of this type.

\item There is a fourth class of representations, more difficult to construct, these are the cuspidal representations. 
We describe here a simple  method of construction of these representations \cite{He}.
The map $N:{\mathbb E}^{\times}\rightarrow {\mathbb F}^{\times}$ induces an injective map $\widehat{{\mathbb F}^{\times}}\hookrightarrow \widehat{{\mathbb E}^{\times}},\mu\mapsto \mu\circ N.$ The elements of $\widehat{{\mathbb E}^{\times}}$ which are not in the image of this map are called primitive characters. An element $\nu$ in $\widehat{{\mathbb E}^{\times}}$  is primitive if and only if $\bar{\nu}=\nu\circ \sigma=\nu^q\not=\nu.$

${\mathbb E}^\times$ acts linearly on ${\mathbb E}$, of  dimension 2 over ${\mathbb F}$, therefore one can identify ${\mathbb E}^\times$ with a subgroup of $GL(2,{\mathbb F}).$ 

Let $Z$ denote the center of $GL(2,{\mathbb F})$ and let $\nu$ be a primitive character of  ${\mathbb E}^\times , $  and $\psi$ a non trivial character of the abelian group $N.$ 
We denote $\nu_\psi$ the character of the subgroup $ZN$ by:
$$\nu_\psi:\begin{pmatrix}a&0\\0&a
\end{pmatrix}n\mapsto \nu(a)\psi(n), a\in {\mathbb F}^\times, n\in N.
$$
The representation $Ind_{ZN}^G(\nu_\psi)$ is, up to equivalence, independent of $\psi$, and we have 
$$Ind_{ZN}^G(\nu_\psi)={\mathcal C}_\nu\oplus Ind_{{\mathbb E}^\times}^G(\nu)$$
where ${\mathcal C}_\nu$ is an irreducible representation  of dimension $q-1,$ the cuspidal representation.

 One has  ${\mathcal C}_{\nu_1}\simeq {\mathcal C}_{\nu_2}$ if and only if $\nu_2=\overline{\nu_1}$  or $\nu_2=\nu_1.$ Therefore $\nu\sim \overline{\nu}.$
There are $\frac{1}{2}q(q-1)$ inequivalent representations of this type.

\end{itemize}

The character table is given by the following table which holds for any $q$ power of prime:
\medskip

\begin{tabular}{|c|c|c|c|c| }\hline
 &$c_1(x) $ &$c_2(x)$ &$c_3(x,y)$&$c_4(\lambda)$\\
\hline\hline
$\chi_\mu$&$\mu(x)^2$&$\mu(x)^2$&$\mu(x)\mu(y)$&$\mu(\lambda\overline{\lambda})$\\
\hline
$I(\mu_1,\mu_2)$&$ (q+1)\mu_1(x)\mu_2(x)$&$\mu_1(x)\mu_2(x)$ &$\mu_1(x)\mu_2(y)+\mu_1(y)\mu_2(x)$&$0$\\
\hline
$St_\mu$&$q\mu(x)^2$&$0$&$\mu(x)\mu(y)$&$-\mu(\lambda\overline{\lambda})$\\
\hline
 ${\mathcal C}_\nu$ &$(q-1)\nu(x)$&$-\nu(x)$&$0$ &$-\nu(\lambda)-\nu(\overline{\lambda})$\\
\hline
\end{tabular}

\subsubsection{Frobenius-Schur indicator}
We give here the list of Frobenius-Schur indicators of irreducible representations of $GL(2,{\mathbb F}).$ Rather than computing these indicators by the formula (\ref{Frobenius}) (which is tedious),  we prefer to use the following description.
We recall that the Frobenius-Schur indicator of a complex  irreducible representation $\pi$ is the integer $\nu_2(\pi)$ which is defined by (\ref{Frobenius}), it takes the value $0$ if the character $\chi_\pi$ is not real, otherwise it takes the value $1$ if the representation is isomorphic to a real representation and is equal to $-1$ if it is not the case.
Therefore by inspection of the character table of $GL(2,{\mathbb F}),$ it is simple to find the representations which have non zero indicator.
We will finally show that there is no irreducible representation having $-1$ value of the  indicator by a counting argument.

\begin{itemize}
\item
$\nu_2(\chi_\mu)$ is non zero if and only if $\mu$ is real. This is possible if and only if $\mu^2=1.$
The solutions of $\mu^2=1$  depend on the parity of $q$. When $q$ is even there is a unique solution $\mu=1=\epsilon_+$, when $q$ is odd one obtains two solutions, the trivial one $\mu=1=\epsilon_+$ and the non trivial character $\epsilon$ of ${\mathbb F}^{\times}$ with value in $\{+1,-1\}.$ We have $\epsilon(x)=1$ if and only if $x$ is a square.

Therefore there are two cases:

if $q$ is even then $\nu_2(\chi_\mu)=0$ if $\mu\not=1$ and $\nu_2(\chi_1)=1$ 

if $q$ is odd then $\nu_2(\chi_\mu)=0$ if $\mu\not\in\{1,\epsilon\}$ and $\nu_2(\chi_1)=\nu_2(\chi_\epsilon)=1.$ 

\item
From the character table, the indicator of  $I(\mu_1,\mu_2)$ is non zero if and only if $\mu_1\mu_2$ and $\mu_1+\mu_2$ are real.
There are two cases:

if $q$ is even then $\nu_2(I(\mu_1,\mu_2))=0$ if and only if $\mu_1\mu_2\not=1.$ 

if $q$ is odd then  $\nu_2(I(\mu_1,\mu_2))=0$  if $\mu_1\mu_2\not=1$  or $\mu_1\mu_2\not=\epsilon.$ If $\mu_1\mu_2=1$ then  $\nu_2(I(\mu_1,\mu_2))\not=0.$
If $\mu_1\mu_2=\epsilon$  then  $\nu_2(I(\mu_1,\mu_2))=0$ unless $\mu_1=1,\mu_2=\epsilon$ or $\mu_2=1,\mu_1=\epsilon.$ We have 
$\nu_2(I(\epsilon,1))=\nu_2(I(1,\epsilon))\in\{-1,+1\}$.

\item
$\nu_2(St_\mu)$ is non zero if and only if $\mu$ is real. This is possible if and only if $\mu^2=1.$
Therefore there are two cases:

if $q$ is even then $\nu_2(St_\mu)=0$ if $\mu\not=1$ and $\nu_2(St)\in\{+1,-1\}$ 

if $q$ is odd then $\nu_2(St_\mu)=0$ if $\mu\not\in\{1,\epsilon\}$ and $\nu_2(St)=\nu_2(St_\epsilon)\in\{+1,-1\}.$ 

\item 
From the character table, the indicator of a cuspidal representation ${\mathcal C}_\rho$ is nul unless $\rho\vert_{{\mathbb F}^\times}^2=1.$
It is easy to see that the condition $\rho(\lambda)+\rho(\overline{\lambda})$ real forces $\rho\vert_{{\mathbb F}^\times}=1,$
and in this case $\rho(\lambda)+\rho(\overline{\lambda})$ is real.
Therefore $\nu_2({\mathcal C}_\rho)=0$ if $\rho\vert_{{\mathbb F}^\times}\not=1,$ and belongs to $\{-1,+1\}$ otherwise.
\end{itemize}
It remains to show that when it is non zero the indicator takes the value $1$.
Let $t$ be the number of non trivial involutions  in $G,$ i.e $t=\vert\{k\in G, k^2=e, k\not=e\}\vert$.
We have $1+t=\sum_{\pi\in \widehat{G}}\nu_2(\pi)dim(\pi)$ using relation (\ref{thetaP2}).
Therefore to show that there is no representation having indicator of value $-1$ it is sufficient to show that the following relation  holds:
\begin{equation}
\label{involution} 1+t=\sum_{\pi\in \widehat{G}, \nu_2(\pi)\not=0}dim(\pi).
\end{equation}
This is done by a counting argument. By solving explicitely  the equation $k^2=e$ in $G$ we obtain that $t=q^2-1$ when $q$ is even and $t=q^2+q+1$ when $q$ is odd. It is easy to compute $\sum_{\pi\in \widehat{G}, \nu_2(\pi)\not=0}dim(\pi)$ from the analysis done above, and we verify that the relation (\ref{involution}) holds.

\subsection{Character table of $PGL(2,{\mathbb F})$}
\subsubsection{Conjugacy classes}

\begin{itemize}
\item
Because $c_1(x)=x I$, these conjugacy classes project on the conjugacy class of $I=\tilde{c}_1$ in $PGL(2,{\mathbb F})$. One conjugacy of this type. 
\item 
Because $c_2(x)$ is similar to $xc_2(1)$ these conjugacy classes project on the conjugacy class of $\tilde{c}_2=\begin{pmatrix}1&1\\0&1
\end{pmatrix}$  in $PGL(2,{\mathbb F})$.One conjugacy class of this type. 

\item
Because $c_3(x,y)=yc_3(x/y,1),$ and $x\not=y$ these conjugacy classes project on the conjugacy classes of $PGL(2,{\mathbb F})$ of the type 
$\begin{pmatrix}x&0\\0&1
\end{pmatrix}=\tilde{c}_3(x)$  with $x\not=1.$ We have $x\backsim x^{-1}.$

Note at this point that there is a distinction between $p$ odd and $p$ even.
When $p$ is $2$ the number of conjugacy classes of this type is $\frac{1}{2}(q-2),$ when $p$ is odd the number of conjugacy classes of this type is $\frac{1}{2}(q-1).$

\item
 $\mu c_4(\lambda)$ is similar to the matrix $c_4(\lambda\mu).$ Therefore to classify the conjugacy classes of $PGL(2,{\mathbb F})$  coming from projection of elements of type $c_4$, it is sufficient to analyze the orbits of ${\mathbb F}^{\times}$ on ${\mathbb E}\setminus {\mathbb F}$ acting by multiplication.
Note also that the ${\mathbb Z}/2{\mathbb Z}$ action $ \lambda\mapsto \overline{\lambda}$ gives the same conjugacy class.
Therefore the set of conjugacy classes of $PGL(2,{\mathbb F})$  coming from projection of elements of type $c_4$, is in bijection with
  (${\mathbb Z}/2{\mathbb Z})\backslash ({\mathbb E}\setminus {\mathbb F})/{\mathbb F}^{\times}.$

 One has to consider the two cases $q$ even and $q$ odd separately.

When $q$ is even, we have $({\mathbb E}\setminus {\mathbb F})/ {\mathbb F}^{\times}\simeq \{\lambda\in {\mathbb E}, \lambda\not=1, N(\lambda)=1\},$ this comes from the fact that ${\mathbb F}\rightarrow {\mathbb F},x\mapsto x^2$ is a bijection in characteristic $2$.
We can always pick a representative of the conjugacy class in $PGL(2,{\mathbb F})$ of 
$\begin{pmatrix}0&-\lambda\overline{\lambda}\\1&\lambda+\overline{\lambda}
\end{pmatrix}$ of the form $\begin{pmatrix}0&1\\1&\lambda+\lambda^{-1}
\end{pmatrix}$  with $\lambda\overline{\lambda}=1.$

Whatever the value of $q$, the $ {\mathbb F}$-linear map  $Tr:{\mathbb E}\rightarrow {\mathbb F}$ is always a  surjection because its kernel  $Ker(Tr)=\{\lambda+\lambda^q=0\}$ has at most $q$ elements.

When $q$ is odd, we can always write $\lambda=x+y\delta$. Therefore in ${\mathbb E}\setminus {\mathbb F}$ there are  two types of orbits. The one with $Tr(\lambda)=0.$ In this case $\lambda=y\delta$, and we can select a unique representative in the ${\mathbb F}^{\times}$  orbit such that $-\lambda\overline{\lambda}=\Delta,$ the representative is therefore 
$\begin{pmatrix}0&\Delta\\1&0
\end{pmatrix}.$ 

In the case of an orbit such that $Tr(\lambda)\not=0,$ we can choose a unique representative in the ${\mathbb F}^{\times}$  orbit  such that $Tr(\lambda)=2.$
We therefore have for this representative  $\lambda=1+y\delta, \overline{\lambda}=1-y\delta, y\not=0.$
The associated element in $PGL(2,{\mathbb F})$ is $\begin{pmatrix}0&-1+y^2\Delta\\1&2
\end{pmatrix}$ which is in the same conjugacy class as $\begin{pmatrix}1&y\Delta\\y&1
\end{pmatrix},$  with $y\not=0$  and $y\backsim y^{-1}.$

\end{itemize}

\subsubsection{Irreducible representations}
Let $Z$ the center of $G$ finite group, a representation $\pi$ of $G$  acting on $V$ is said to admit a central character $\rho:Z\rightarrow {\mathbb C}^{\times}$ if $\pi(z)=\rho(z)Id_V,\forall z\in Z.$ If $\pi$ admits a central character, we will denote ${\mathfrak Z}_\pi$ its central character.
From Schur lemma, finite dimensional irreducible representations of $G$ admit central character. Moreover if $\pi$ and $\pi'$  have central character, $\pi\otimes \pi'$ admits central character ${\mathfrak Z}_{\pi\otimes \pi'}={\mathfrak Z}_\pi {\mathfrak Z}_{\pi'},$ and $\pi\oplus\pi'$ admits a central character if and only if ${\mathfrak Z}_\pi ={\mathfrak Z}_{\pi'}.$ 

There is a bijection between representions of $PGL(2,{\mathbb F})$ and representations of $GL(2,{\mathbb F})$ of trivial central character.
Therefore the  finite dimensional irreducible representations of $PGL(2,{\mathbb F})$ are in bijection with the finite dimensional irreducible representations of $GL(2,{\mathbb F})$ with trivial center character.
Therefore the classification of  representations of $PGL(2,{\mathbb F})$ is  obtained from the character table of $GL(2,{\mathbb F})$ by  selecting the irreducible representations of  $GL(2,{\mathbb F})$ where the value of $c_1(x)$ is trivial.

One dimensional representations of $PGL(2,{\mathbb F})$ are of the type $\chi_\mu$ with $\mu^2=1.$

The irreducible representations $PGL(2,{\mathbb F})$  belong to the following list:

\begin{itemize}
\item $\chi_\mu$ with $\mu^2=1$
\item $I(\mu_1,\mu_2)$ with $\mu_1\mu_2=1.$
\item $St_\mu$ with $\mu^2=1.$
\item ${\mathcal C}_\nu$ with $\nu\vert_{{\mathbb F}^{\times}}=id$
\end{itemize}

When $p$ is 2, the list of inequivalent representations of $PGL(2,{\mathbb F})\simeq SL(2,{\mathbb F})$ is:
\begin{itemize}
\item $\chi_1$,  the trivial representation of dimension 1.
\item $I(\mu,1)$ of dimension $q+1$ with $\mu\not=1$, $\mu\sim\mu^{-1}$. There are $\frac{1}{2}q-1$ inequivalent representations of this type.
\item $St$ of dimension $q$.
\item ${\mathcal C}_\nu$ of dimension $q-1$ with $\nu\vert_{{\mathbb F}^{\times}}=id$ and $\nu\sim\nu^{-1}.$ There are $\frac{1}{2}q$ inequivalent representations of this type.
\end{itemize}

When $p$ is odd , the list of inequivalent representations of  $PGL(2,{\mathbb F})$ is:
\begin{itemize}
\item $\chi_1, \chi_{\epsilon}$  inequivalent  representation of dimension 1.
\item $I(\mu,1)$ of dimension $q+1$ with $\mu\not\in \{1,\epsilon\}$, $\mu\sim\mu^{-1}$. There are $\frac{1}{2}(q-3)$ inequivalent representations of this type.
\item $St,St_{\epsilon} =St\otimes \chi_{\epsilon}$ inequivalent representations of dimension $q$.
\item ${\mathcal C}_\nu$ of dimension $q-1$ with $\nu\vert_{{\mathbb F}^{\times}}=id$ and $\nu\sim\nu^{-1}.$ There are $\frac{1}{2}(q-1)$ inequivalent representations of this type.
\end{itemize}

Here are the character table of $PGL(2,{\mathbb F})$ in the two cases:

Case where  $p=2:$

\begin{tabular}{|c|c|c|c|c| }\hline
 &$\begin{pmatrix}1&0\\0&1
\end{pmatrix}$  &  $\begin{pmatrix}1&1\\0&1
\end{pmatrix}$  &$\begin{pmatrix}x&0\\0&1
\end{pmatrix},x\backsim x^{-1}$ &  $\begin{pmatrix}0&1\\1&\lambda+\overline{\lambda}
\end{pmatrix},\overline{\lambda}=\lambda^{-1}\backsim \lambda$ \\
\hline\hline
$\chi_1$&$1$&$1$&$1$&$1$\\
\hline
$I(\mu,1)$&$ (q+1)$&$1$ &$\mu(x)+\mu(x^{-1})$&$0$\\
\hline
$St$ &$q$&$0$&$1$&$-1$\\
\hline
 ${\mathcal C}_\nu$ &$(q-1)$&$-1$&$0$ &$-\nu(\lambda)-\nu(\overline{\lambda})$\\
\hline
\end{tabular}
\bigskip

Case where $p$ is odd:

\begin{tabular}{|c|c|c|c|c| }\hline
 &$\begin{pmatrix}1&0\\0&1
\end{pmatrix}$  &  $\begin{pmatrix}1&1\\0&1
\end{pmatrix}$  &$\begin{pmatrix}x&0\\0&1
\end{pmatrix},x\backsim x^{-1}$ &  $\begin{pmatrix}0&-\lambda\overline{\lambda}\\1&\lambda+\overline{\lambda}
\end{pmatrix},\overline{\lambda}=\lambda^{-1}\backsim \lambda$ \\
\hline\hline
$\chi_1$&$1$&$1$&$1$&$1$\\
\hline
$\chi_{\epsilon}$&$1$&$1$&$\epsilon(x)$&$\epsilon(\lambda\overline{\lambda})$\\
\hline
$I(\mu,1)$&$ (q+1)$&$1$ &$\mu(x)+\mu(x^{-1})$&$0$\\
\hline
$St$ &$q$&$0$&$1$&$-1$\\
\hline
$St_{\epsilon}$ &$q$&$0$&$\epsilon(x)$&$-\epsilon(\lambda\overline{\lambda})$\\
\hline
 ${\mathcal C}_\nu$ &$(q-1)$&$-1$&$0$ &$-\nu(\lambda)-\nu(\overline{\lambda})$\\
\hline
\end{tabular}
\bigskip

\subsubsection{Frobenius-Schur indicator}
There is a bijection between the irreducible representations of $PGL(2,{\mathbb F})$ and the irreducible representations of $GL(2,{\mathbb F})$ having trivial central character. moreover it is easy to show that if $\tilde{\pi}$ is the irreducible representation of  $PGL(2,{\mathbb F})$ associated, via this bijection , to the irreducible representation $\pi$ of $GL(2,{\mathbb F})$, we have $\nu_2(\pi)=\nu_2(\tilde{\pi}).$
From the explicit knowledge of the Frobenius-Schur indicator in the $GL(2,{\mathbb F})$ case, one obtains the value of the Frobenius-Schur indicator in the $PGL(2,{\mathbb F})$ case.
We easily obtain that  every irreducible representation of $PGL(2,{\mathbb F})$ has  Schur-Frobenius indicator equal to $1.$

\subsection{Decomposition of tensor products of irreducible representations of  of $GL(2,{\mathbb F})$}
We answer here the problem of computing  explicitely the multiplicities of irreducible representations in the tensor product of two irreducible representation of $G=GL(2,{\mathbb F})$. This  problem is equivalent to the explicit description of the fusion ring of 
$GL(2,{\mathbb F}).$ This is a simple computation, involving only sum over roots of unity, but the explicit result seems not to be known by the experts. 
Because irreducible representations of $PGL(2,{\mathbb F})$ are in bijection with the subset of  irreducible representations of $GL(2,{\mathbb F})$ of trivial central character,
the fusion ring of $PGL(2,{\mathbb F})$ is a subring of the the fusion ring of $GL(2,{\mathbb F})$. We therefore study only the structure of the fusion ring of $GL(2,{\mathbb F}).$

If $\pi$ is an irreducible representation, let $\check{\pi}$ be the contragredient representation. 
Let $\pi,\pi', \pi"$ irreducible finite dimensional complex representations of $G$, we denote 
\begin{eqnarray}&&\langle\pi\pi'\rangle=\frac{1}{\vert G\vert}\sum_{g\in G}\chi_{\pi}(g)\chi_{\pi'}(g),\\
&&\langle\pi\pi'\pi''\rangle=\frac{1}{\vert G\vert}\sum_{g\in G}\chi_{\pi}(g)\chi_{\pi'}(g)\chi_{\pi''}(g).
\end{eqnarray}

The multiplicities are non negative integers which can be computed by the formula 
\begin{eqnarray}
N_{\pi,\pi'}^{\pi''}=\frac{1}{\vert G\vert}\sum_{g\in G}\chi_{\pi}(g)\chi_{\pi'}(g)\chi_{\pi''}(g^{-1}).
\end{eqnarray}
We prefer  to compute the numbers $N_{\pi,\pi',\pi''}=N_{\pi,\pi'}^{\check{\pi}''}= \langle\pi\pi'\pi''\rangle$  which are symmetric under the exchange of $\pi,\pi',\pi''$  and for which the explicit formula are simple.

We have 
\begin{equation}
\check{\chi}_\mu\simeq\chi_{\mu^{-1}},
\check{I}(\mu_1,\mu_2)\simeq I(\mu_1^{-1},\mu_2^{-1}),
\check{St}_\mu\simeq St_{\mu^{-1}},
\check{\mathcal C}_\nu\simeq {\mathcal C}_{\nu^{-1}}.
\end{equation}

We will denote  $\delta$ the delta function on $\widehat{{\mathbb F}^{\times}} $ located on the unit element whereas  $\delta^{{\mathbb E}^{\times }}$ denotes the delta function on $\widehat{{\mathbb E}^{\times}} $ located on the unit element. We will denote if $\mu\in\widehat{{\mathbb F}^{\times}} $ and $\nu\in \widehat{{\mathbb E}^{\times}} $, for typographic reasons, $\delta_\mu=\delta(\mu), \delta_\nu^{{\mathbb E}^{\times }}=\delta^{{\mathbb E}^{\times }}(\nu).$
If $\nu$ is a character of ${\mathbb E}^{\times }$ we will also denote $\underline{\nu}=\nu\vert_{{\mathbb F}^{\times} }\in \widehat{{\mathbb F}^{\times}}.$

We give here the list of non zero elements of $\langle\pi\pi'\rangle$ for $\pi,\pi'$ irreducible representations of $G,$

\begin{proposition}

\begin{eqnarray*}
&&\langle \chi_\mu\chi_\nu \rangle=\delta_{\mu\nu},\\
&&\langle I(\mu_1,\mu_2)  I(\nu_1,\nu_2)    \rangle =\delta_{\mu_1\nu_1}\delta_{\mu_2\nu_2}+\delta_{\mu_1\nu_2}\delta_{\mu_2\nu_1},\\
&&\langle St_\mu St_\nu\rangle=\delta_{\mu\nu},\\
&&\langle {\mathcal  C}_\nu {\mathcal  C}_\rho\rangle=\delta_{\nu\rho}+\delta_{\nu\overline{\rho}}.
\end{eqnarray*}
\end{proposition}

Proof: it follows from the orthogonality relations on characters and of the stated above equivalence of representations.

\begin{proposition}

\begin{eqnarray*}
&&\langle I(\mu_1,\mu_2)I(\nu_1,\nu_2)I(\rho_1,\rho_2)\rangle=\\
&&\delta_{\mu_1\nu_1\rho_1\mu_2\nu_2\rho_2}+
\delta_{\mu_1\nu_1\rho_1}\delta_{\mu_2\nu_2\rho_2}+
\delta_{\mu_2\nu_1\rho_1}\delta_{\mu_1\nu_2\rho_2}+\delta_{\mu_1\nu_2\rho_1}\delta_{\mu_2\nu_1\rho_2}+\delta_{\mu_1\nu_1\rho_2}\delta_{\mu_2\nu_2\rho_1},\\
&&\langle I(\mu_1,\mu_2)I(\nu_1,\nu_2)St_\rho\rangle=\delta_{\mu_1\mu_2\nu_1\nu_2\rho^2}+\delta_{\mu_1\nu_1\rho}\delta_{\mu_2\nu_2\rho}+
\delta_{\mu_1\nu_2\rho}\delta_{\mu_2\nu_1\rho},\\
&&\langle  I(\mu_1,\mu_2)I(\nu_1,\nu_2){\mathcal C}_\rho\rangle=\delta_{\mu_1\mu_2\nu_1\nu_2\underline{\rho}},\\
&&\langle I(\mu_1,\mu_2) St_\nu St_\rho\rangle=\delta_{\mu_1\mu_2\nu^2\rho^2}+\delta_{\mu_1\nu\rho}\delta_{\mu_2\nu\rho},\\
&&\langle St_{\mu} St_{\nu} St_{\rho}\rangle=\delta_{\mu^2\nu^2\rho^2},\\
&&\langle I(\mu_1,\mu_2) St_\nu {\mathcal C}_\rho\rangle=\delta_{\mu_1\mu_2\nu^2\underline{\rho}},\\
&&\langle I(\mu_1,\mu_2) {\mathcal C}_\nu    {\mathcal C}_\rho  \rangle=\delta_{\mu_1\mu_2\underline{\nu}\underline{\rho}},\\
&&\langle St_{\mu} St_{\nu}{\mathcal C}_\rho\rangle= \delta_{\mu^2\nu^2\underline{\rho}},\\
&&\langle  St_{\mu} {\mathcal C}_\nu  {\mathcal C}_\rho\rangle=\delta_{\mu^2\underline{\nu}\underline{\rho}}-\delta^{{\mathbb E}^{\times }}_{\mu \overline {\mu}{\nu}{\rho}}-\delta^{{\mathbb E}^{\times }}_{\mu \overline {\mu}{\nu}\overline{{\rho}}},\\
&&\langle   {\mathcal C}_\mu   {\mathcal C}_\nu  {\mathcal C}_\rho\rangle=
\delta_{\underline{\mu}\underline{\nu}\underline{\rho}}-(\delta^{{\mathbb E}^{\times }}_{\mu\nu\rho}+\delta^{{\mathbb E}^{\times }}_{\overline{\mu}\nu\rho}+\delta^{{\mathbb E}^{\times }}_{\mu\overline{\nu}\rho}+\delta^{{\mathbb E}^{\times }}_{\mu\nu\overline{\rho}}).
\end{eqnarray*}
(These last  two expressions on the right, although conveniently written with a minus sign, are non negative integer belonging to \{0,1\}).
\end{proposition}

Proof:

This follows from  direct computations.
We only give the proof of the first identity, the others follow from the same type of proof.
\begin{eqnarray*}
&&\vert G\vert \langle I(\mu_1,\mu_2)I(\nu_1,\nu_2)I(\rho_1,\rho_2)\rangle=\\
&&=\sum_{x\in  {\mathbb F}^{\times}}(q+1)^3(\mu_1\mu_2\nu_1\nu_2\rho_1\rho_2)(x) +
\sum_{x\in  {\mathbb F}^{\times}}(\mu_1\mu_2\nu_1\nu_2\rho_1\rho_2)(x) (q^2-1)+q(q+1)\times\\
&&\sum_{\substack{\{x,y\}\subset  {\mathbb F}^{\times}\\x\not=y }} (\mu_1(x)\mu_2(y)+\mu_1(y)\mu_2(x))(\nu_1(x)\nu_2(y)+\nu_1(y)\nu_2(x))(\rho_1(x)\rho_2(y)+\rho_1(y)\rho_2(x)).
\end{eqnarray*}
expanding all the product and using the identity
\begin{equation}
\sum_{x\in \mathbb{F}^{\times}}\mu(x)=(q-1)\delta_\mu,  \mu\in\widehat{{\mathbb F}^{\times}},
\end{equation}
 we obtain the announced result.
$\Box$

This is the complete list of $N_{\pi,\pi',\pi''}$ where none of the $\pi,\pi',\pi''$ are one dimensional representations.
Using 
\begin{eqnarray*}
&&\chi_\mu\otimes\chi_\nu\simeq\chi_{\mu\nu},\;\;
\chi_\mu\otimes I(\nu,\rho)\simeq I(\mu\nu, \mu\rho)\\
&&\chi_\mu\otimes St_\nu\simeq St_{\mu\nu},\;\;
\chi_\mu\otimes {\mathcal C}_\nu\simeq {\mathcal C}_{\mu\nu},
\end{eqnarray*}

we obtain that if   $\pi''$ is one dimensional,  we have $ N_{\pi,\pi',\pi''}= \langle \pi (\pi'\otimes\pi'')\rangle $ already computed.

\end{appendix}

\end{document}